\documentclass{amsart}

\usepackage{macros}
\standardmargins\standardpackages\standardrenews\standardlabeling
\colorcommentstrue

\draftfalse

\begin{document}
\title{An analogue of a theorem of Kurzweil}

\authordavid

\begin{arxivabstract}
A theorem of Kurzweil ('55) on inhomogeneous Diophantine approximation states that if $\theta$ is an irrational number, then the following are equivalent: (A) for every decreasing positive function $\psi$ such that $\sum_{q = 1}^\infty \psi(q) = \infty$, and for almost every $s\in\mathbb R$, there exist infinitely many $q\in\mathbb N$ such that $\|q\theta - s\| < \psi(q)$, and (B) $\theta$ is badly approximable. This theorem is not true if one adds to condition (A) the hypothesis that the function $q\mapsto q\psi(q)$ is decreasing. In this paper we find a condition on the continued fraction expansion of $\theta$ which is equivalent to the modified version of condition (A). This expands on a recent paper of D. H. Kim ('14).
\end{arxivabstract}

\maketitle

An irrational number $\theta$ is said to be \emph{badly approximable} (or \emph{of bounded type}) if there exists $\epsilon > 0$ such that for every rational $p/q\in\Q$,
\[
\left|\theta - \frac{p}{q}\right| \geq \frac{\epsilon}{q^2}\cdot
\]
It is well-known that an irrational number $\theta$ is badly approximable if and only if the partial quotients of $\theta$ form a bounded sequence. Another equivalent condition was given by Kurzweil \cite{Kurzweil2}. To state it, let us define the set
\[
W(\theta,\psi) = \{s\in\R : \exists^\infty q\in\N \; \|q\theta - s\| < \psi(q)\},
\]
where $\|\cdot\|$ denotes distance to the nearest integer. Then Kurzweil's result may be stated as follows: $\theta$ is badly approximable if and only if for every decreasing function $\psi:\N\to(0,\infty)$ such that $\sum_{q = 1}^\infty \psi(q) = \infty$, the set $W(\theta,\psi)$ has full measure. (Note that if $\sum_{q = 1}^\infty \psi(q) < \infty$, then the set $W(\theta,\psi)$ has measure zero by the Borel--Cantelli lemma.)

Rather than considering all decreasing functions $\psi$, one may consider the smaller class of \emph{Khinchin sequences}: a function $\psi:\N\to(0,\infty)$ is called a Khinchin sequence if, in addition to the divergence condition $\sum_{q = 1}^\infty \psi(q) = \infty$, the function $q\mapsto q\psi(q)$ is nonincreasing. Although less natural than the condition that $\psi$ is decreasing, the hypothesis that a sequence is a Khinchin sequence is significant both for historical reasons (Khinchin first proved his eponymous theorem \cite{Khinchin1} in the setting of Khinchin sequences, although his theorem was later generalized) and because such sequences are often easier to work with.

Let $\theta$ be an irrational number and let $\psi$ be a Khinchin sequence. A recent paper of D. H. Kim \cite{Kim2} gives a criterion, based on the continued fraction expansion of $\theta$, for the set $W(\theta,\psi)$ to have full measure.\arxivfootnote{After this paper was written, Kim extended his result to all positive decreasing sequences in a joint paper with M. Fuchs \cite{FuchsKim}.} However, his paper leaves open the question of finding an analogue of Kurzweil's theorem in the setting of Khinchin sequences, although he proves several results in that direction \cite[\63]{Kim2}. The aim of this paper is to complete the work of Kim by proving such an analogue.

{\bf Acknowledgements.} The author thanks Lior Fishman for helpful comments.

\section{Statement of results}
We first recall the main theorem of \cite{Kim2}, rephrased slightly:\arxivfootnote{Technically, the result of \cite{Kim2} applies to the sets $\bigcap_{\epsilon > 0} W(\theta,\epsilon\psi)$ and not directly to the sets $W(\theta,\psi)$. But since the convergence or divergence of the series \eqref{kim} is invariant under a slight perturbation of $\psi$, \cite[Theorem 2.1]{Kim2} and Theorem \ref{theoremkim} are equivalent.}
\begin{theorem}[{\cite[Theorem 2.1]{Kim2}}]
\label{theoremkim}
Fix $\theta\in\R\butnot\Q$ and let $(q_k)_0^\infty$ be the sequence of the denominators of the convergents of $\theta$. Let $\psi:\N\to (0,\infty)$ be a Khinchin sequence, and let $\phi(q) = 1/(q\psi(q))$. Then the following are equivalent:
\begin{itemize}
\item[(A)] $W(\theta,\psi)$ has full measure.
\item[(B)] The series
\begin{equation}
\label{kim}
\sum_{k = 0}^\infty \frac{\log\phi(q_k)\wedge\log(q_{k + 1}/q_k)}{\phi(q_k)}
\end{equation}
diverges. (In this paper, $\wedge$ and $\vee$ denote minimum and maximum, respectively.)
\end{itemize}
\end{theorem}

To state our main theorem, we use the notation
\[
\Sigma\big((a_i)_1^n : m\big)
\]
to denote the sum of the $m$ largest elements of the sequence $(a_i)_1^n$, with $\Sigma\big((a_i)_1^n : m\big) = \sum_1^n a_i$ if $m\geq n$. For $\alpha\geq 0$, we let $\Sigma\big((a_i)_1^n : \alpha\big) = \Sigma\big((a_i)_1^n : \lfloor\alpha\rfloor\big)$.

\begin{theorem}
\label{maintheorem}
Fix $\theta\in\R\butnot\Q$ and let $(q_k)_0^\infty$ be the sequence of the denominators of the convergents of $\theta$. Then the following are equivalent:
\begin{itemize}
\item[(A)] For every Khinchin sequence $\psi:\N\to (0,\infty)$, the set $W(\theta,\psi)$ has full measure.
\item[(B)] For some $\epsilon > 0$,
\[
\limsup_{k\to\infty} \frac{1}{\log(q_k)} \Sigma\left(\left(\log\left(\frac{q_{i + 1}}{q_i}\right)\right)_{i = 0}^{k - 1} : \frac{\epsilon\log(q_k)}{\log\log(q_k)}\right) < 1.
\]
\end{itemize}
\end{theorem}

\begin{remark*}
Since condition (B) of Theorem \ref{maintheorem} is not equivalent to the condition that the sequence $(q_k)_1^\infty$, it follows from Kurzweil's theorem that condition (A) is not equivalent to the condition that for every decreasing positive function $\psi:\N\to (0,\infty)$ such that $\sum_{q = 1}^\infty \psi(q) = \infty$, the set $W(\theta,\psi)$ has full measure. In particular, there exists a decreasing positive function $\psi:\N\to(0,\infty)$ such that $\sum_{q = 1}^\infty \psi(q) = \infty$ and such that there is no Khinchin sequence $\psi':\N\to(0,\infty)$ with $\psi'(q)\leq \psi(q)$ for all $q$. An example of such a sequence is given by the formula
\[
\psi(q) = 2^{-n_k} \;\; (2^{n_{k - 1}} \leq q < 2^{n_k})
\]
where $(n_k)_1^\infty$ is any sequence of integers such that $n_k - n_{k - 1} \geq k$ for all $k$.
\end{remark*}

\section{Proof of Theorem \ref{maintheorem}}
{\bf Convention.} The symbol $\asymp$ will denote a coarse multiplicative asymptotic, i.e. $A_n\asymp B_n$ means that there exists a constant $C > 0$ (the \emph{implied constant}) such that $C^{-1} B_n \leq A_n\leq C B_n$.

\begin{proof}[Proof of \text{(A) \implies (B)}]
By contradiction, suppose that (B) is false. Then for each $n\in\N$, there exists $k_n\in\N$ such that
\[
\frac{1}{\log(q_{k_n})} \Sigma\left(\left(\log\left(\frac{q_{i + 1}}{q_i}\right)\right)_{i = 0}^{k_n - 1} : \frac1{2^n} \frac{\log(q_{k_n})}{\log\log(q_{k_n})}\right) \geq 1 - \frac1{2^n}\cdot
\]
Without loss of generality, suppose that the sequence $(k_n)_1^\infty$ is increasing, and let $k_0 = 0$. For each $n\geq 1$, let $S_n'$ be a subset of $\{0,\ldots,k_n - 1\}$ of cardinality at most $\frac1{2^n} \frac{\log(q_{k_n})}{\log\log(q_{k_n})}$ such that
\[
\sum_{k\in S_n'} \log\left(\frac{q_{k + 1}}{q_k}\right) \geq (1 - 2^{-n})\log(q_{k_n}).
\]
Then let $S_n = S_n'\butnot\{0,\ldots,k_{n - 1} - 1\}$ and $T_n = \{k_{n - 1},\ldots,k_n - 1\}\butnot S_n$. Then
\begin{equation}
\label{cardSn}
\#(S_n) \leq \frac1{2^n} \frac{\log(q_{k_n})}{\log\log(q_{k_n})}
\end{equation}
and
\begin{equation}
\label{massTn}
\sum_{k\in T_n} \log\left(\frac{q_{k + 1}}{q_k}\right) \leq \log(q_{k_n}) - \sum_{k\in S_n'} \log\left(\frac{q_{k + 1}}{q_k}\right) \leq 2^{-n}\log(q_{k_n}).
\end{equation}
Now define the function $\phi:\N\to(0,\infty)$ by the formula
\[
\phi(q) = \log(q_{k_n}) \all q_{k_{n - 1}} \leq q < q_{k_n}.
\]
Then $\phi$ is nondecreasing, and
\begin{align*}
\sum_{q = 1}^\infty \frac{1}{q\phi(q)}
&= \sum_{n = 1}^\infty \frac{1}{\log(q_{k_n})}\sum_{q = q_{k_{n - 1}}}^{q_{k_n} - 1} \frac{1}{q}
\asymp \sum_{n = 1}^\infty \frac{\log(q_{k_n}/q_{k_{n - 1}})}{\log(q_{k_n})}\\
&= \sum_{n = 1}^\infty \left[1 - \frac{\log(q_{k_{n - 1}})}{\log(q_{k_n})}\right]
\asymp \sum_{n = 1}^\infty 1\wedge\log\left(\frac{\log(q_{k_n})}{\log(q_{k_{n - 1}})}\right) = \infty.
\end{align*}
Thus $\psi(q) = 1/(q\phi(q))$ is a Khinchin sequence. So by (A) together with Theorem \ref{theoremkim}, the series \eqref{kim} diverges. On the contrary, we show that \eqref{kim} converges:
\begin{align*}
&\sum_{k = 0}^\infty \frac{\log\phi(q_k)\wedge\log(q_{k + 1}/q_k)}{\phi(q_k)}\\
&\leq \sum_{n = 1}^\infty \left[\sum_{k\in S_n} \frac{\log\phi(q_k)}{\phi(q_k)} + \sum_{k\in T_n} \frac{\log(q_{k + 1}/q_k)}{\phi(q_k)}\right] \noreason\\
&= \sum_{n = 1}^\infty \left[\frac{\log\log(q_{k_n})}{\log(q_{k_n})}\#(S_n) + \frac{1}{\log(q_{k_n})}\sum_{k\in T_n} \log(q_{k + 1}/q_k)\right] \noreason\\
&\leq \sum_{n = 1}^\infty \left[\frac{1}{2^n} + \frac{1}{2^n}\right] \by{\eqref{cardSn} and \eqref{massTn}}\\
&= 2 < \infty.
\end{align*}
This contradiction completes the proof.
\end{proof}

\begin{proof}[Proof of \text{(B) \implies (A)}]
Let $\psi:\N\to\infty$ be a Khinchin sequence, and by contradiction suppose that $W(\theta,\psi)$ does not have full measure. Then by Theorem \ref{theoremkim}, the series \eqref{kim} converges, where $\phi(q) = 1/(q\psi(q))$ is nondecreasing. Let
\[
S = \{k : \phi(q_k) \leq q_{k + 1}/q_k\}, \; T = \N\butnot S,
\]
so that
\[
\infty > \sum_{k = 0}^\infty \frac{\log\phi(q_k)\wedge\log(q_{k + 1}/q_k)}{\phi(q_k)} = \sum_{k\in S} \frac{\log\phi(q_k)}{\phi(q_k)} + \sum_{k\in T} \frac{\log(q_{k + 1}/q_k)}{\phi(q_k)}\cdot
\]
For each $m\in\N$, let $Q_m$ be the largest integer such that $\phi(Q_m)\leq 2^m$. Then
\begin{align*}
\frac{1}{\phi(q)} &\asymp \sum_{\substack{m\in\N : 2^m \geq \phi(q)}} \frac{1}{2^m} = \sum_{\substack{m\in\N : q\leq Q_m}} \frac{1}{2^m}\\
\frac{\log\phi(q)}{\phi(q)} &\asymp \sum_{\substack{m\in\N : 2^m \geq \phi(q)}} \frac{m}{2^m} = \sum_{\substack{m\in\N : q\leq Q_m}} \frac{m}{2^m}
\end{align*}
and thus
\begin{align*}
\infty &= \sum_{q = 1}^\infty \frac{1}{q\phi(q)} \asymp \sum_{m = 0}^\infty \sum_{q = 1}^{Q_m} \frac{1}{q 2^m} \asymp \sum_{m = 0}^\infty \frac{\log(Q_m)}{2^m}\\
\infty &> \sum_{k\in S} \frac{\log\phi(q_k)}{\phi(q_k)} + \sum_{k\in T} \frac{\log(q_{k + 1}/q_k)}{\phi(q_k)}\\
&\asymp \sum_{m = 0}^\infty \left[\frac{m}{2^m}\#\{k\in S: q_k\leq Q_m\} + \frac{1}{2^m} \sum_{\substack{k\in T \\ q_k\leq Q_m}}\log\left(\frac{q_{k + 1}}{q_k}\right)\right].
\end{align*}
It follows that if
\begin{align*}
\lambda_m &= \frac{\displaystyle\frac{m}{2^m}\#\{k\in S: q_k\leq Q_m\} + \frac{1}{2^m} \sum_{\substack{k\in T \\ q_k\leq Q_m}}\log\left(\frac{q_{k + 1}}{q_k}\right)}{\displaystyle\frac{\log(Q_m)}{2^m}}\\
&= \frac{1}{\log(Q_m)}\left[m\#\{k\in S: q_k\leq Q_m\} + \sum_{\substack{k\in T \\ q_k\leq Q_m}}\log\left(\frac{q_{k + 1}}{q_k}\right)\right],
\end{align*}
then
\[
\liminf_{m\to\infty} \lambda_m = 0.
\]
On the other hand, if
\[
\kappa_m = \frac{m}{2^m}\#\{k\in S: q_k\leq Q_m\},
\]
then
\[
\lim_{m\to\infty}\kappa_m = 0.
\]
Fix $\epsilon > 0$, and choose $m\geq 2$ such that $\lambda_m , \kappa_m \leq \epsilon$. Then
\[
\frac{m}{2^m}\vee\frac{m}{\log(Q_m)} \leq \frac{\epsilon}{\#\{k\in S: q_k\leq Q_m\}}\cdot
\]
Consider the function
\[
f(x) = \frac{x}{2^x}\vee\frac{x}{\log(Q_m)} \hspace{.2 in} (x\geq 2).
\]
Since $f$ is the maximum of an increasing function and a decreasing function, $f$ has a unique minimum, which occurs when the two inputs to the maximum agree, namely at $x = \log_2\log(Q_m)$. Thus
\[
\frac{\epsilon}{\#\{k\in S: q_k\leq Q_m\}} \geq f(m) \geq \min(f) = \frac{\log_2\log(Q_m)}{\log(Q_m)} \geq \frac{\log\log(Q_m)}{\log(Q_m)}
\]
i.e.
\[
\#\{k\in S: q_k\leq Q_m\} \leq \frac{\epsilon\log(Q_m)}{\log\log(Q_m)}\cdot
\]
On the other hand, since $\lambda_m\leq\epsilon$,
\[
\sum_{\substack{k\in T \\ q_k\leq Q_m}}\log\left(\frac{q_{k + 1}}{q_k}\right) \leq \epsilon\log(Q_m).
\]
Let $k_m$ be the smallest integer such that $Q_m < q_{k_m}$. Then
\begin{align*}
\#\{k\in S: k < k_m\} &\leq \frac{\epsilon\log(q_{k_m})}{\log\log(q_{k_m})}\\
\sum_{\substack{k\in T \\ k < k_m}}\log\left(\frac{q_{k + 1}}{q_k}\right) &\leq \epsilon\log(q_{k_m})
\end{align*}
and thus
\[
\Sigma\left(\left(\log\left(\frac{q_{i + 1}}{q_i}\right)\right)_{i = 0}^{k_m - 1} : \frac{\epsilon\log(q_{k_m})}{\log\log(q_{k_m})}\right) \geq (1 - \epsilon)\log(q_{k_m}).
\]
Since $\epsilon$ was arbitrary and $k_m\to\infty$, for all $\epsilon > 0$ we have
\[
\limsup_{k\to\infty} \frac{1}{\log(q_k)}\Sigma\left(\left(\log\left(\frac{q_{i + 1}}{q_i}\right)\right)_{i = 0}^{k - 1} : \frac{\epsilon\log(q_k)}{\log\log(q_k)}\right) = 1,
\]
contradicting (B).
\end{proof}

\section{Consequences of Theorem \ref{maintheorem}}
In this section we use Theorem \ref{maintheorem} to prove some necessary and sufficient conditions on $\theta$ for $W(\theta,\psi)$ to be full measure for every Khinchin sequence $\psi$, including reproving some results from \cite[\63]{Kim2}. For convenience let
\[
\Omega = \{\theta\in\R : \text{for every Khinchin sequence $\psi$, the set $W(\theta,\psi)$ has full measure}\}.
\]
In other words, $\Omega$ is the set of all $\theta$ such that the equivalent conditions of Theorem \ref{maintheorem} hold.
\begin{theorem}
Fix $\theta\in\R\butnot\Q$ and let $(q_k)_0^\infty$ be the sequence of the denominators of the convergents of $\theta$.
\begin{itemize}
\item[(i)] If
\begin{equation}
\label{cond1}
\limsup_{k\to\infty} \frac{\log(q_k)}{k} < \infty,
\end{equation}
then $\theta\in\Omega$.
\item[(ii)] If
\begin{equation}
\label{cond2}
\limsup_{k\to\infty} \frac{\log(q_k)}{k\log(k)} = \infty,
\end{equation}
then $\theta\notin\Omega$.
\item[(iii)] If
\begin{equation}
\label{cond3}
\sum_{k = 2}^\infty \frac{1}{\log(q_k)} < \infty
\end{equation}
then $\theta\notin\Omega$.
\item[(iv)] If
\begin{equation}
\label{cond4}
\limsup_{k\to\infty}\frac{q_{k + 1}/q_k}{\log(q_k)} < \infty,
\end{equation}
then $\theta\in\Omega$.
\item[(v)] If
\begin{equation}
\label{cond5}
\limsup_{k\to\infty}\frac{\log(q_{k + 1}/q_k)}{\log(q_k)} = \infty,
\end{equation}
then $\theta\notin\Omega$.
\end{itemize}
\end{theorem}
\begin{remark*}
Parts (i), (iii), and (iv) correspond to \cite[Theorem 3.1 and Proposition 3.2]{Kim2}. Although in some cases the new proofs are not shorter than the old proofs, having two proofs may bring further insight.
\end{remark*}
\begin{remark*}
By well-known facts about continued fractions (e.g. \cite[Theorems 9 and 13]{Khinchin_book}), the conditions \eqref{cond4} and \eqref{cond5} have interpretations in terms of Diophantine approximation:
\begin{itemize}
\item $\theta$ satisfies \eqref{cond4} if and only if for some $\epsilon > 0$, $\theta$ is not $\psi$-approximable, where
\[
\psi(q) = \frac{\epsilon}{q^2\log(q)}\cdot
\]
We recall that a number $\theta$ is called \emph{$\psi$-approximable} if there exist infinitely many rationals $p/q\in\Q$ such that
\[
\left|\theta - \frac pq\right| < \psi(q).
\]
\item $\theta$ satisfies \eqref{cond5} if and only if $\theta$ is a Liouville number. We recall that a number $\theta$ is called \emph{Liouville} if for all $n\in\N$, $\theta$ is $\psi_n$-approximable, where $\psi_n(q) = q^{-n}$.
\end{itemize}
\end{remark*}
\begin{remark*}
Any badly approximable number $\theta$ satisfies both \eqref{cond1} and \eqref{cond4}, so $\BA\subset\Omega$. This can also be seen from Kurzweil's theorem.
\end{remark*}
\begin{remark*}
The continued fraction expansion of $e$ (see e.g. \cite{Cohn}) satisfies \eqref{cond4}, so $e\in\Omega$.
\end{remark*}
\begin{proof}[Proof of \text{(i)}]
Choose $M < \infty$ so that for all $k$, $\log(q_k) \leq Mk$. Let $\epsilon > 0$ be arbitrary (e.g. $\epsilon = 1$). Then for sufficiently large $k$,
\[
\frac{\epsilon\log(q_k)}{\log\log(q_k)} \leq \frac{\epsilon Mk}{\log(Mk)} \leq \frac{k}{8}\cdot
\]
Let $S \subset \{0,\ldots,k - 1\}$ be a subset of cardinality at most $k/8$, and let $T = \{0,\ldots,k - 1\}\butnot S$. A counting argument shows that
\[
\#\{i = 0,\ldots,k - 1 \text{ even} : i, i + 1\in T\} \geq k/4,
\]
and thus
\[
\sum_{i\in T} \log(q_{i + 1}/q_i) \geq \sum_{\substack{i \text{ even} \\ i,i + 1\in T}} \log(q_{i + 2}/q_i) \geq (k/4) \log(2) \geq \frac{\log(2)}{4M}\log(q_k).
\]
It follows that
\[
\frac{1}{\log(q_k)} \Sigma\left(\left(\log\left(\frac{q_{i + 1}}{q_i}\right)\right)_{i = 0}^{k - 1} : \frac{\epsilon\log(q_k)}{\log\log(q_k)}\right) \leq 1 - \frac{\log(2)}{4M}\cdot
\]
To complete the proof, we take the limsup as $k\to\infty$ and then apply Theorem \ref{maintheorem}.
\end{proof}
\begin{proof}[Proof of \text{(ii)}]
Fix $\epsilon > 0$. By assumption, there exist infinitely many $k$ satisfying
\[
\log(q_k) \geq \frac{2}{\epsilon}k\log(k).
\]
For such $k$,
\[
\frac{\epsilon\log(q_k)}{\log\log(q_k)} \geq \frac{\epsilon(2/\epsilon)k\log(k)}{\log\big((2/\epsilon)k\log(k)\big)} \geq \frac{2 k\log(k)}{\log(k^2)} = k,
\]
where the middle inequality holds for all $k$ sufficiently large. But then
\[
\Sigma\left(\left(\log\left(\frac{q_{i + 1}}{q_i}\right)\right)_{i = 0}^{k - 1} : \frac{\epsilon\log(q_k)}{\log\log(q_k)}\right) = \sum_{i = 0}^{k - 1} \log\left(\frac{q_{i + 1}}{q_i}\right) = \log(q_k).
\]
To complete the proof, we divide by $\log(q_k)$, take the limsup as $k\to\infty$, and apply Theorem \ref{maintheorem}.
\end{proof}
Since \eqref{cond3} implies \eqref{cond2}, (iii) does not require a separate proof.
\begin{proof}[Proof of \text{(iv)}]
Choose $M < \infty$ such that for all $k$, $q_{k + 1}/q_k \leq M\log(q_k)$. Then for all $\epsilon > 0$ and $k\in\N$,
\begin{align*}
\Sigma\left(\left(\log\left(\frac{q_{i + 1}}{q_i}\right)\right)_{i = 0}^{k - 1} : \frac{\epsilon\log(q_k)}{\log\log(q_k)}\right)
&\leq \frac{\epsilon\log(q_k)}{\log\log(q_k)} \max\left\{\log\left(\frac{q_{i + 1}}{q_i}\right) : i = 0,\ldots,k - 1\right\}\\
&\leq \frac{\epsilon\log(q_k)}{\log\log(q_k)} \log(M\log(q_k)) \leq 2\epsilon\log(q_k),
\end{align*}
where the last inequality holds for all $k$ large enough such that $q_k\geq e^M$. To complete the proof, we let $\epsilon = 1/4$, divide by $\log(q_k)$, take the limsup as $k\to\infty$, and apply Theorem \ref{maintheorem}.
\end{proof}
\begin{proof}[Proof of \text{(v)}]
The assumption \eqref{cond5} implies that
\[
\limsup_{k\to\infty} \frac{\log(q_k/q_{k - 1})}{\log(q_k)} = 1.
\]
Fix $\epsilon > 0$. By assumption, there exist infinitely many $k$ such that 
\[
\log(q_k/q_{k - 1}) \geq (1 - \epsilon)\log(q_k).
\]
For such $k$, if we assume that $k$ is chosen large enough so that $\frac{\epsilon\log(q_k)}{\log\log(q_k)}\geq 1$, then
\[
\Sigma\left(\left(\log\left(\frac{q_{i + 1}}{q_i}\right)\right)_{i = 0}^{k - 1} : \frac{\epsilon\log(q_k)}{\log\log(q_k)}\right)
\geq \log\left(\frac{q_k}{q_{k - 1}}\right) \geq (1 - \epsilon)\log(q_k).
\]
To complete the proof, we divide by $\log(q_k)$, take the limsup as $k\to\infty$, use the fact that $\epsilon$ was arbitrary, and apply Theorem \ref{maintheorem}.
\end{proof}

\bibliographystyle{amsplain}

\bibliography{bibliography}

\end{document}